\documentclass{article}%
\usepackage{amsmath}
\usepackage{amsfonts}
\usepackage{amssymb}
\usepackage{amsthm}
\usepackage{fullpage}
\usepackage{graphicx}
\usepackage{color}
\usepackage{authblk}
\usepackage{tikz}
\usepackage{comment}
\usepackage{subcaption}%
\usepackage{hyperref}
\providecommand{\U}[1]{\protect\rule{.1in}{.1in}}
\providecommand{\U}[1]{\protect\rule{.1in}{.1in}}
\allowdisplaybreaks
\newtheorem{theorem}{Theorem}

\newtheorem{corollary}[theorem]{Corollary}

\newtheorem{lemma}[theorem]{Lemma}
\newtheorem{notation}[theorem]{Notation}
\newtheorem{proposition}[theorem]{Proposition}

\theoremstyle{definition}
\newtheorem{definition}[theorem]{Definition}
\newtheorem{example}[theorem]{Example}

\newtheorem{remark}[theorem]{Remark}

\begin{document}

\title{Stratifying Discriminant Hypersurface}
\author{\mbox{}\\\vspace{-1em} Rizeng Chen \\School of Mathematical Sciences, Peking University, Beijing 100091, China \\xiaxueqaq@stu.pku.edu.cn \\\vspace{1em} Hoon Hong \\Department of Mathematics, North Carolina State University, Raleigh NC 27695,
USA \\hong@ncsu.edu \\\vspace{1em} Jing Yang\thanks{Corresponding author.} \\SMS--School of Mathematical Sciences, Guangxi Minzu University, Nanning
530006, China \\yangjing0930@gmail.com }
\date{}
\maketitle

\begin{abstract}
\noindent This paper investigates the stratification of the discriminant
hypersurface associated with a univariate polynomial via the number of its
distinct complex roots. We introduce two novel approaches different from the
one based on subdiscriminants. The first approach stratifies the discriminant
hypersurface by recursively removing all the \emph{lowest-order} points, while
the second one stratifies the discriminant hypersurface by recursively
removing all the \emph{smooth} points. Both approaches rely solely on the
discriminant itself instead of using high-order subdiscriminants. These
results offer new insights into the intrinsic geometry of the discriminant and
its connection to root multiplicity.

\end{abstract}

\section{Introduction}

Discriminant of a univariate polynomial is a fundamental object in
computational algebraic geometry. It plays an essential role on understanding
the root structure of a given polynomial. In the past centuries, people
proposed various versions of discriminants, extended to subdiscriminants, and
developed rich theories on them~(to list a few, \cite{1853_Sylvester,
1864_Sylvester, 1864_Sylvester_extension, 1967_Collins,
1994_Gelfand_Kapranov_Zelevinsky, 2009_Lazard_McCallum, 2009_Buse_Mourrain,
2016_Buse_Karasoulou, 2017_Buse_Nonkane, 2013_Esterov,
2023_Dickenstein_DiRocco_Morrison, 1996_Yang_Hou_Zeng, 1999_Yang,
2021_Hong_Yang, 2024_Hong_Yang:non-nested, 2012_Nie} ). Then they applied the
theories to tackle various fundamental problems in computational algebraic
geometry (see, e.g., \cite{2002_Xia_Yang, 2001_Yang_Hou_Xia, 2005_Yang_Xia,
2007_Lazard_Rouillier, 2012_Maza_Xia_Xiao}), such as finding condition on the
coefficients of the univariate polynomial so that the polynomial has a
specified number of distinct complex/real roots, or has a specified
multiplicities, and so on (see, e.g., \cite{1853_Sylvester,
1996_Yang_Hou_Zeng, 1998_Gonzalez_Vega_Recio_Lombardi_Roy, 2021_Hong_Yang,
2024_Hong_Yang:non-nested}).

We revisit the following classical problem: \emph{Stratify}  the discriminant hypersurface associated with a univariate polynomial via the number of its distinct complex roots.   It is well known that such stratification can be obtained via  subdiscriminants. In this paper, we introduce two novel approaches different from the one based on subdiscriminants.

To illustrate the stratification problem and our contributions, consider the
monic quartic polynomial $F=x^{4}+a_{3}x^{3}+a_{2}x^{2}+a_{1}x+a_{0}$
where~$a_{i}$'s are indeterminates/parameters. Its discriminant $D$ is given
by
\begin{align*}
D  &  =256a_{0}^{3}-27a_{1}^{4}-128a_{0}^{2}a_{2}^{2}-192a_{0}^{2}a_{1}%
a_{3}+144a_{0}a_{1}^{2}a_{2}\\
&  -4a_{1}^{2}a_{2}^{3}+18a_{1}^{3}a_{2}a_{3}-6a_{0}a_{1}^{2}a_{3}^{2}\\
&  +16a_{0}a_{2}^{4}+144a_{0}^{2}a_{2}a_{3}^{2}-80a_{0}a_{1}a_{2}^{2}%
a_{3}+18a_{0}a_{1}a_{2}a_{3}^{3}\\
&  -4a_{1}^{3}a_{3}^{3}-4a_{0}a_{2}^{3}a_{3}^{2}-27a_{0}^{2}a_{3}^{4}%
+a_{1}^{2}a_{2}^{2}a_{3}^{2}.
\end{align*}
Define
\begin{equation}
Z_{k}=\left\{  (a_{0},a_{1},a_{2},a_{3})\in\mathbb{C}^{4}:x^{4}+a_{3}%
x^{3}+a_{2}x^{2}+a_{1}x+a_{0}=0 \text{ has at most }k \text{ distinct
roots}\right\}  .
\end{equation}
It is obvious that $Z_{4}=\mathbb{C}^{4}$ and~$Z_{0}=\varnothing$. Also, it is
well-known that $Z_{k}$'s are algebraic sets defined by the subdiscriminants.
Moreover, the zero locus of $D$ is the hypersurface $Z_{3}$ in the parameter
space. It is easy to construct monic quartics having 1, 2, 3, 4 distinct roots
respectively:%
\[%
\begin{array}
[c]{llll}%
x^{4}-4x^{3}+6x^{2}-4x+1 & = & \left(  x-1\right)  ^{4} & \in Z_{1}\backslash
Z_{0}\\
x^{4}-5x^{3}+9x^{2}-7x+2 & = & \left(  x-1\right)  ^{3}\left(  x-2\right)  &
\in Z_{2}\backslash Z_{1}\\
x^{4}-7x^{3}+17x^{2}-17x+6 & = & \left(  x-1\right)  ^{2}\left(  x-2\right)
\left(  x-3\right)  & \in Z_{3}\backslash Z_{2}\\
x^{4}-10x^{3}+35x^{2}-50x+24 & = & \left(  x-1\right)  \left(  x-2\right)
\left(  x-3\right)  \left(  x-4\right)  & \in Z_{4}\backslash Z_{3}%
\end{array}
\]
Therefore, we conclude that
\[
\ \varnothing=Z_{0}\subsetneq Z_{1}\subsetneq Z_{2}\subsetneq Z_{3}=V(D) \subsetneq
Z_{4}=\mathbb{C}^{4}.
\]
In general, let $F=x^{n}+a_{n-1}x^{n-1}+\cdots+a_{0}$ be a general monic
univariate polynomial of degree $n$ and let $D$ be the discriminant of $F$. We
have
\[
\varnothing=Z_{0}\subsetneq Z_{1}\subsetneq Z_{2}\subsetneq\cdots\subsetneq Z_{n-1}%
=V(D)\subsetneq Z_{n}=\mathbb{C}^{n}.
\]
which can be illustrated by the following diagram.

\begin{center}
\begin{tikzpicture}
\draw[draw=black]        (-2.5,-2.0) rectangle  ++(5.0,4.6);
\draw[draw=black]        (0, 0.1) circle [radius=2.0];
\draw[draw=black]        (0,-0.1) circle [radius=1.6];
\draw[draw=black]        (0,-0.3) circle [radius=1.2];
\draw[draw=black,dotted] (0,-0.5) circle [radius=0.8];
\draw[draw=black]        (0,-0.7) circle [radius=0.4];
\node (C4) at ( 0, 2.30){$Z_n=\mathbb{C}^n$};
\node (V4) at ( 0, 1.80){$Z_{n-1}$};
\node (V3) at ( 0, 1.20){$Z_{n-2}$};
\node (V2) at ( 0, 0.60){$Z_{n-3}$};
\node (V1) at ( 0,-0.55){$Z_{1}$};
\node[circle, fill=black, minimum size=3pt, inner sep=0pt, outer sep=0pt]
(d4) at (1.3, 2.00){};%
\node[circle, fill=black, minimum size=3pt, inner sep=0pt, outer sep=0pt]
(d3) at (1.0, 1.50){};%
\node[circle, fill=black, minimum size=3pt, inner sep=0pt, outer sep=0pt]
(d2) at (0.7, 0.95){};
\node[circle, fill=black, minimum size=3pt, inner sep=0pt, outer sep=0pt]
(d1) at (0.1,-0.75){};
\node (w4)     at (3.00, 2.00)[right]{$n\ \ \ \ \ $ distinct roots};
\node (w3)     at (3.00, 1.50)[right]{$n-1$ distinct roots};
\node (w2)     at (3.00, 0.95)[right]{$n-2$ distinct roots};
\node (wdots)  at (3.00, 0.25)[right]{$\vdots$};
\node (w1)     at (3.00,-0.75)[right]{$1\ \ \ \ \ $ distinct root};
\draw[thin] (d4) -- (w4);
\draw[thin] (d3) -- (w3);
\draw[thin] (d2) -- (w2);
\draw[thin] (d1) -- (w1);
\end{tikzpicture}

\end{center}

\noindent Note that $Z_{0},\ldots,Z_{n}$ provides a stratification of
$\mathbb{C}^{n}$ such that, in each strata $S_{i}=Z_{i}-Z_{i-1},$ the number
of distinct roots of $F$ is invariant, namely $i$. Naturally we would like to
describe the stratification in terms of the coefficients of the polynomial
$F$, that is, to find defining polynomials of $Z_{i}$ in terms of the coefficients.

Note that $Z_{n}=\mathbb{C}^{n}$ and $Z_{0}=\varnothing$ are trivial. What is
really interesting is the remaining sets $Z_{1},\ldots,Z_{n-1}$. It is
well-known that $Z_{i}$ can be described via sub-discriminants (see Section
\ref{sec:review}). 

For a long time, researchers have believed that the discriminant itself should carry the same amount of information as subdiscriminants. For example, Collins used all subdiscriminants as part of projection polynomials in cylindrical algebraic decomposition but conjectured that, in most cases, only the discriminant is necessary \cite{collins1975quantifier}. This insight motivated McCallum to refine projection techniques, successfully designing a much smaller set of projection polynomials that rely on discriminants rather than all subdiscriminants \cite{mccallum1998improved}. Yet, despite these advances, it remains unclear how to determine the number of distinct roots using only the discriminant.

This naturally leads to the core question underlying our work: is it possible to describe $Z_{i}$ in
terms of the coefficients without using sub-discriminants and instead using
only \emph{discriminant itself}? If so, how? The main contribution of this
paper is to provide two different approaches to solving this problem.

\begin{description}
\item[Contribution 1.] We partition the space of the coefficients into subsets
according to the order of discriminant on the points, so that the order of the
discriminant is invariant over each subset. All the elements in each subset
have the same order and generate polynomials with the same number of distinct
roots. 
In summary, this approach stratifies the discriminant hypersurface by recursively removing all the \emph{lowest-order} points.
See Theorem \ref{thm:order}.

\item[Contribution 2.] We view the hypersurface defined by the discriminant as
a variety in $\mathbb{C}^{n}$ and analyze the singular locus of this variety,
the singular locus of the singular variety obtained in the previous step, ...,
and so on. This yields a stratification of the discriminant hypersurface,
whose strata are smooth varieties. We prove that each stratum exactly
corresponds to the set of polynomials having a given number of distinct roots.
In summary, this approach stratifies the discriminant
hypersurface by recursively removing all the \emph{smooth} points. See
Theorem \ref{thm:singular_stratification}.
\end{description}

\begin{figure}[ptbh]
\centering
\begin{subfigure}[t]{0.38\textwidth}
\centering
\includegraphics[width=\linewidth]{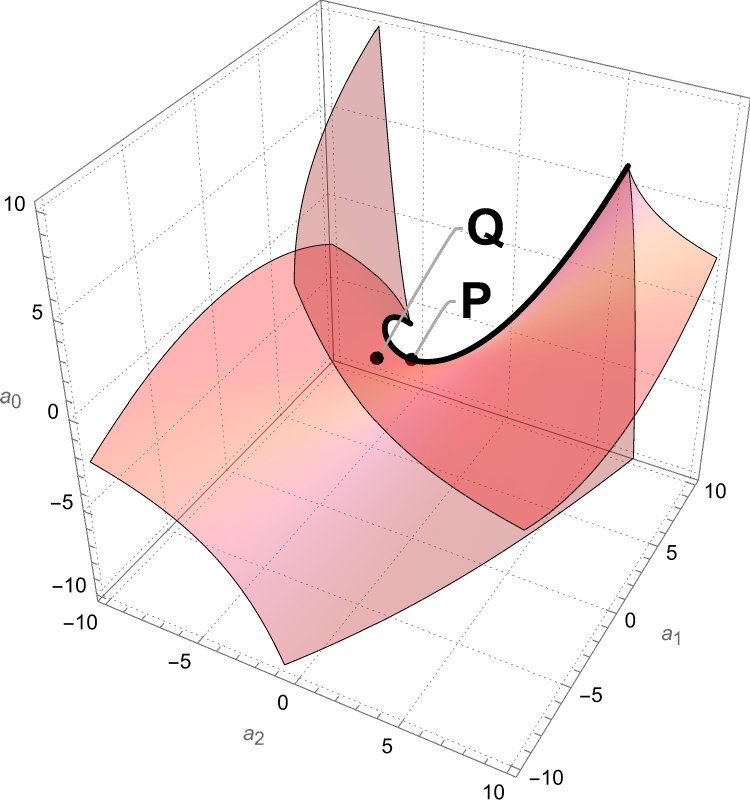}
\caption{ The discriminant hypersurface $Z_2$}
\end{subfigure}
\begin{subfigure}[t]{0.6\textwidth}
\centering
\begin{tikzpicture}
\draw[draw=black] (-2.15,-2.0) rectangle  ++(3.6,4.1);
\draw[draw=black] (-0.3,-0.2) circle [radius=1.6];
\draw[draw=black] (-0.3,-0.8) circle [radius=0.9];
\node (C3) at (-0.2, 1.8){$Z_3=\mathbb{C}^3$};
\node (V2) at (-0.2, 1.13){$Z_2$};
\node (V1) at (-0.2, -0.18){$Z_1$};
\node[circle, fill=black, minimum size=3pt, inner sep=0pt, outer sep=0pt]
(d3) at (1.2,1.2){};
\node[circle, fill=black, minimum size=3pt, inner sep=0pt, outer sep=0pt]
(d2) at (0.7,0.2){};
\node[circle, fill=black, minimum size=3pt, inner sep=0pt, outer sep=0pt]
(d1) at (0.2,-0.8){};
\node (w3) at (2.1, 1.2)[right]{3 distinct roots};
\node (w2) at (2.1, 0.2)[right]{2 distinct roots};
\node (w1) at (2.1,-0.8)[right]{1 distinct roots};
\draw[thin] (d3) -- (w3);
\draw[thin] (d2) -- (w2);
\draw[thin] (d1) -- (w1);
\draw[draw=black] (-2.2,-2.2);
\end{tikzpicture}
\caption{A diagram showing the case $n=3$}%
\label{fig:diagram-n-equal-to-3}%
\end{subfigure}
\caption{Stratification for $n=3$}%
\label{fig:monic-cubic}%
\end{figure}
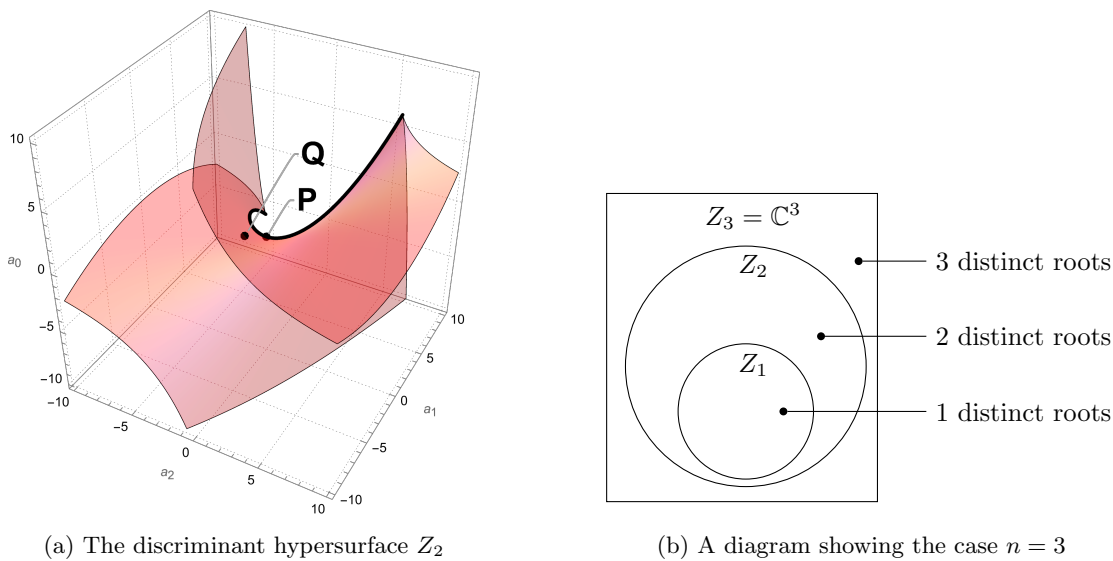

\begin{figure}[ptbh]
\centering
\begin{subfigure}[t]{0.38\textwidth}
\centering
\includegraphics[width=\linewidth]{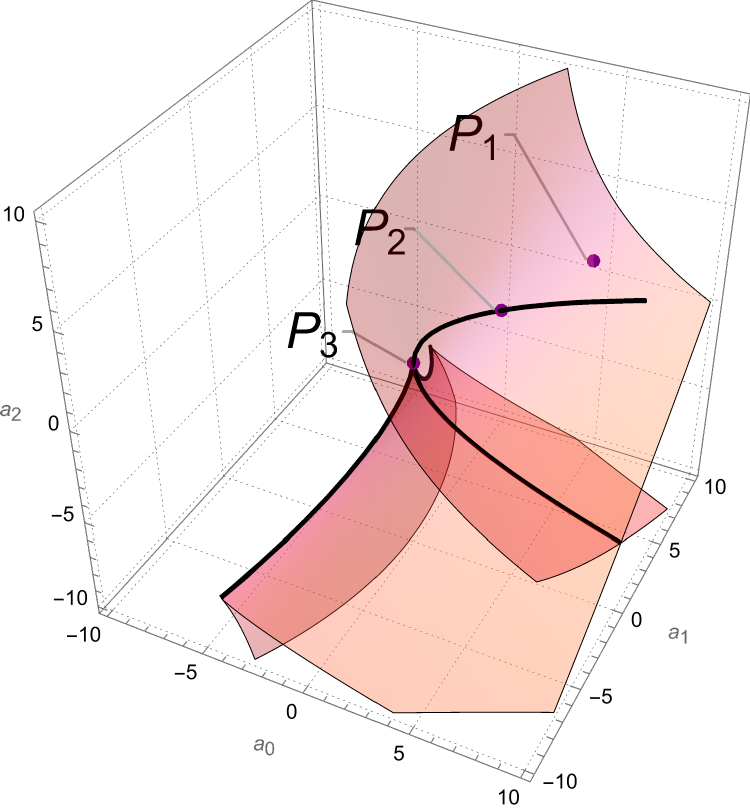}
\caption{A cross section of the discriminant hypersurface $Z_3$}
\label{fig:a3-is-zero}
\end{subfigure}
\begin{subfigure}[t]{0.6\textwidth}
\centering
\begin{tikzpicture}
\draw[draw=black] (-2.12,-2) rectangle  ++(3.6,4.1);
\draw[draw=black] (-0.3,-0.2) circle [radius=1.6];
\draw[draw=black] (-0.3,-0.7) circle [radius=1.0];
\draw[draw=black] (-0.3,-1.1) circle [radius=0.5];
\node (C4) at (-0.3,1.80) {$Z_4=\mathbb{C}^4$};
\node (V3) at (-0.3,1.10) {$Z_3$};
\node (V2) at (-0.3,0.00) {$Z_2$};
\node (V1) at (-0.3,-0.85){$Z_1$};
\node[circle, fill=black, minimum size=3pt, inner sep=0pt, outer sep=0pt] (d4) at (1.2,1.2){};
\node[circle, fill=black, minimum size=3pt, inner sep=0pt, outer sep=0pt] (d3) at (0.7,0.4){};
\node[circle, fill=black, minimum size=3pt, inner sep=0pt, outer sep=0pt] (d2) at (0.2,-0.4){};
\node[circle, fill=black, minimum size=3pt, inner sep=0pt, outer sep=0pt] (d1) at (-0.3,-1.2){};
\node (w4) at (2.1,1.2)[right]{4 distinct roots};
\node (w3) at (2.1,+0.4)[right]{3 distinct roots};
\node (w2) at (2.1,-0.4)[right]{2 distinct roots};
\node (w1) at (2.1,-1.2)[right]{1 distinct root};
\draw[thin] (d4) -- (w4);
\draw[thin] (d3) -- (w3);
\draw[thin] (d2) -- (w2);
\draw[thin] (d1) -- (w1);
\end{tikzpicture}
\caption{A diagram showing the case $n=4$}%
\label{fig:diagram}%
\end{subfigure}
\caption{Stratification for $n=4$}%
\label{fig:monic-quartic}%
\end{figure}
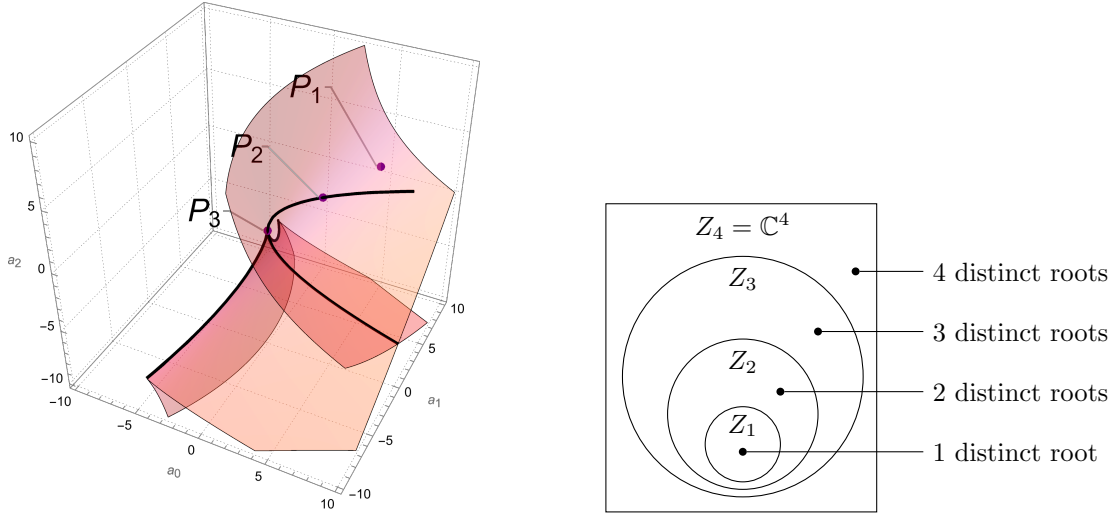

\begin{example}
To better illustrate our contribution, we will begin with a simple example,
where
\[
F=x^{3}+a_{2}x^{2}+a_{1}x+a_{0}.
\]
The discriminant is
\[
D=-4a_{2}^{3}a_{0}+a_{2}^{2}a_{1}^{2}+18a_{2}a_{1}a_{0}-4a_{1}^{3}-27a_{0}%
^{2}.
\]
We plot its zero locus in Figure \ref{fig:monic-cubic}(a) in pink. As one may
tell from the figure, the zero locus of $D$ is a surface in the
three-dimensional space, and it is singular along the black curve. We can
arbitrarily pick one singular point and one smooth point from $V(D)$, say
$P=(0,0,0)$ and $Q=(0,-3,2)$, corresponding to $x^{3}$ and $x^{3}%
-3x+2=(x-1)^{2}(x+2)$ respectively.

It is easy to verify that $D$ is of order $2$ at $P$, and of order $1$ at $Q$.
In fact,
\[%
\begin{array}
[c]{rcl}%
D & = & -27 a_{0}^{2}+\left(  -4 a_{1}^{3}+18 a_{2} a_{1} a_{0}+a_{2}^{2}
a_{1}^{2}-4 a_{2}^{3} a_{0}\right) \\
& = & -108 a+(\text{terms of total degree}\geq2\text{ in }a, b+3, c-2).
\end{array}
\]

\noindent This agrees with what is predicted in our results. On the one hand,
the polynomial $x^{3}$ has only one root counted without multiplicity, so the
order of $D$ is $2$ at $P$ and $P$ is a singularity of $V(D)$. On the other
hand, the polynomial $x^{3}-3x+2$ has two distinct roots, so the order of $D$
is $1$ at $Q$ and $Q$ is a smooth point of~$V(D)$.
\end{example}

\begin{example}
\label{ex:quartic} Let us continue with the monic quartic example. Since it is
impossible to directly visualize the 4-dimensional space, we will plot a
section of the discriminant hypersurface $Z_{3}$, cut by the hyperplane
$a_{3}=0$. This is still quite informative, because one can always apply the
Tschirnhausen transformation to eliminate the $a_{3}x^{3}$ term. See
Figure~\ref{fig:a3-is-zero}.

The points corresponding to
\[%
\begin{array}
[c]{lclcl}%
P_{1} & = & x^{4}+2x^{2}+8 x+5 & = & (x+1)^{2}(x-1+2\mathrm{i}%
)(x-1-2\mathrm{i}),\\
P_{2} & = & x^{4}+4 x^{2}+4 & = & (x-\sqrt{2}\mathrm{i})^{2}(x+\sqrt
{2}\mathrm{i})^{2},\\
P_{3} & = & x^{4} &  &
\end{array}
\]
are also plotted in Figure \ref{fig:a3-is-zero}. One can see that while
$P_{1}$ is lying on the pink surface, one of the black curves, which are
``self-intersections'' of the pink surface, passes through $P_{2}$. Also,
$P_{3}$ is the intersection of the two black curves. Moreover, it can be
easily checked that the discriminant $D$ is of order $1$ at $P_{1}$, order $2
$ at $P_{2}$ and order $3$ at $P_{3}$.

In fact, this is not just a coincidence. Since $P_{1}$ has three distinct
roots, our Theorem \ref{thm:order} shows that $D$ is of order $1$ at $P_{1}$,
and Theorem \ref{thm:singular_stratification} concludes that $P_{1}$ is a
smooth point of $Z_{3}$. Similarly, $P_{2}$ has only two distinct roots, so
$D$ is of order $2$ at it and it must be a singular point of $Z_{3}$, lying on
the black curves. As for the last point, $D$ is of order $3$ at $P_{3}$ and it
is a singularity of the singular locus $\mathop{\mathrm{Sing}}Z_{3}$, with $\operatorname{Sing} V$ represents the singular locus of the variety $V$ hereinafter.

\end{example}

\begin{remark}
\ 

\begin{itemize}
\item The combination of our two results together shows that the smooth
stratification of the discriminant hypersurface can be directly read off from
the information of order.

\item While the smooth stratification of the discriminant hypersurface is an
intrinsic object, the order of the discriminant polynomial depends on the
embedding in the ambient space. One may wonder whether our first contribution
can be formulated intrinsically. The answer is affirmative because the notion
of order can be reinterpreted using commutative algebra. More explicitly, the
order of a polynomial $F$ at a point $p\in V(F)$ coincides with the
Hilbert-Samuel multiplicity of the local ring $\mathcal{O}_{V(F),p}$. See, for
example, \cite[Corollary 5.5.9]{2008_Greuel_Gerhard}.

In this paper, we will however stick to the classical definition of order
because it is conceptually simpler, hence more accessible to the general readers.
\end{itemize}

\end{remark}

\paragraph{Related Work:}

\begin{enumerate}
\item Let $F\left(  a,x\right)  ~$and $G\left(  b,x\right)  $ be two
univariate polynomials in $x$ with indeterminate coefficients $a$ and $b$
respectively. Let $R\left(  a,b\right)  $ be the resultant of $F\left(
a,x\right)  $ and $G\left(  b,x\right)  $ with respect to $x$. Let $\alpha
,\beta$ be numeric coefficient vectors such that $\deg_{x}F\left(
\alpha,x\right)  =\deg_{x}F\left(  a,x\right)  $ and $\deg_{x}G\left(
\beta,x\right)  =\deg_{x}G\left(  b,x\right)  $. Recall the following
well-known facts:

\begin{enumerate}
\item The number of distinct roots of $F\left(  \alpha,x\right)  $ is
$n-\deg\gcd\left(  F\left(  \alpha,x\right)  ,F^{\prime}\left(  \alpha,x\right)
\right)  $.

\item The order of $R$ at $\left(  \alpha,\beta\right)  $ is equal to $\deg
\gcd\left(  F\left(  \alpha,x\right)  ,G\left(  \beta,x\right)  \right)  $. In
\cite[Appendix]{2018_Hong_Sendra}, Bus\'{e} provided a short and elegant proof.

\item The discriminant of $F\left(  \alpha,x\right)  $ is the resultant of
$F\left(  \alpha,x\right)  $ and $F^{\prime}\left(  \alpha,x\right)  $ up to a
non-zero constant.
\end{enumerate}

Thus one might suspect that our Contribution 1 could be easily obtained from
specializing the work of Bus\'{e}. However, it is not clear whether this is
possible because 
\begin{itemize}
\item order is not preserved under specialization, and
\item the partial derivatives  (with respect to coefficients)  of resultant are different from the partial derivatives of discriminant.
\end{itemize}

\item Let $Y_{r}(m,n)$ be the determinantal variety consisting of all $m\times
n$ matrices of rank $\leq r$. Recall the following classical results:

\begin{enumerate}
\item For every $1\leq r\leq l=\min(m,n)-1,$ the variety $Y_{r-1}(m,n)$
coincides with the singular locus of~$Y_{r}(m,n)$ \cite[Proposition 3.2,
Chapter 4]{1994_Gelfand_Kapranov_Zelevinsky}. Hence,
\[
Y_{l}(m,n)\supsetneq\cdots\supsetneq Y_{1}(m,n)\supsetneq Y_{0}(m,n)
\]
yields a smooth stratification of $Y_{l}(m,n)$.

\item The number of distinct roots of $F$ is uniquely determined by the rank
of the Sylvester matrix of~$F\ $and $F^{\prime}$.
\end{enumerate}

Thus one may wonder whether our Contribution 2 could be easily obtained from
specializing the smooth stratification of $Y_{l}\left(  m,n\right)  $.
However, it is not clear whether this is possible because singularity is not
preserved under specialization.

\item The notion of iterated discriminants was considered by several authors
\cite{2023_Dickenstein_DiRocco_Morrison,2016_Han,2009_Buse_Mourrain,2009_Lazard_McCallum,1868_Henrici}%
. It provides some information about the singular locus of the discriminant
hypersurface. For example, \cite[Proposition 9]{2009_Lazard_McCallum}
characterizes $\mathop{\mathrm{Sing}}Z_{n-1}$. However it is not clear how to
obtain our Contribution 2 because their result do not show how to characterize
the other $Z_{i}$'s.

\item Let $C_{\mu}$ be the family of polynomials of type $\left(
x-x_{1}\right)  ^{\mu_{1}}\cdots\left(  x-x_{k}\right)  ^{\mu_{k}}$ where
$\mu=\left(  \mu_{1},\ldots,\mu_{k}\right)  $ is a partition of $n$ where
again $x_{i}$'s are not necessarily distinct. It is a folklore that the
irreducible components of $Z_{k}$ consist of $C_{\mu}$'s where $\#\mu=k$ (see
our Corollary \ref{cor:irr-dec-of-v-k}). Therefore researches on $Z_{k}$ or
$C_{\mu}$ are closely intertwined.

The object~$C_{\mu}$ have been studied under several distinct terminologies:
coincident root locus \cite{2003_Chipalkatti,2012_Kurmann}, multiple root
locus \cite{2016_Lee_Sturmfels}, $\lambda$-Chow variety \cite{2012_oeding},
etc. In \cite{2021_Hong_Yang,2024_Hong_Yang:non-nested}, Hong and Yang
discovered \textquotedblleft simple\textquotedblright\ defining polynomials
for $C_{\mu}$'s under the name \textquotedblleft multiplicity
discriminants\textquotedblright. In~\cite{2003_Chipalkatti}, Chipalkatti fully
determined the singular locus of $C_{\mu}$. Later, the object
$\mathop{\mathrm{Sing}}C_{\mu}$ was given another equivalent characterization
by Kurmann in \cite{2012_Kurmann}.

Based on Chipalkatti's results and a careful analysis on all the
partitions~$\mu$ by their largest number, we managed to prove our Contribution 2.
\end{enumerate}

The current paper is structured as follows. In Section \ref{sec:sdisc}, we give a brief review on the stratification of dicriminant hypersurface based on subdiscriminants. In Section \ref{sec:order}, we
prove that the discriminant hypersurface can be stratified based on the
order of the discriminant at the points so that the order of the discriminant
remains constant within each stratum. In Section \ref{sec:smooth}, we present another
stratification of the discriminant hypersurface based on smoothness where each stratum exactly corresponds to the set of polynomials
having a given number of distinct roots. The paper is concluded in Section
\ref{sec:conclusion} with some further remarks.

We assume that the reader is familiar with the following notations: order of
polynomial at a point and singular locus of a variety.

\section{Subdiscriminant-based Stratification}\label{sec:sdisc}

\label{sec:review} \label{sec:sdis}

We briefly review the classical notion of subdiscriminants and
their  use in stratification~\cite{1853_Sylvester,1864_Sylvester}.

\begin{definition}
[Subdiscrimiant]Let $F\left(  x\right)  =\sum\limits_{i}a_{i}x^{i}$ be a monic
polynomial of degree $n$, that is $a_{n}=1$. Let~$F^{\prime}\left(  x\right)
=\sum\limits_{i}b_{i}x^{i}$ be the derivative of $F$, that is $b_{i}=\left(
i+1\right)  a_{i+1}$. Then the $k$-th subdiscriminant of $F$, written as
$D_{k}$, is defined by
\[
D_{k}=\left(  -1\right)  ^{\binom{n-k}{2}}\det\left[
\begin{array}
[c]{cccccc}%
a_{n} & a_{n-1} & \cdots & \cdots & \cdots & \cdots\\
& \ddots & \ddots &  &  & \\
&  & a_{n} & a_{n-1} & \cdots & \cdots\\
b_{n-1} & b_{n-2} & \cdots & \cdots & \cdots & \cdots\\
& \ddots & \ddots &  &  & \\
&  & \ddots & \ddots &  & \\
&  &  & b_{n-1} & b_{n-2} & \cdots
\end{array}
\right]
\]
where \# of the top rows (involving $a$'s) is $n-1-k$ and the \# of the bottom
rows (involving $b$'s) is $n-k$.
\end{definition}

\begin{example}
Let $n=4$. Recall that $a_{4}=1$. We have\
\begin{align*}
D_{0}  &  =\left(  -1\right)  ^{\frac{4\cdot3}{2}}\left\vert
\begin{array}
[c]{ccccccc}%
a_{4} & a_{3} & a_{2} & a_{1} & a_{0} &  & \\
& a_{4} & a_{3} & a_{2} & a_{1} & a_{0} & \\
&  & a_{4} & a_{3} & a_{2} & a_{1} & a_{0}\\
4a_{4} & 3a_{3} & 2a_{2} & 1a_{1} &  &  & \\
& 4a_{4} & 3a_{3} & 2a_{2} & 1a_{1} &  & \\
&  & 4a_{4} & 3a_{3} & 2a_{2} & 1a_{1} & \\
&  &  & 4a_{4} & 3a_{3} & 2a_{2} & 1a_{1}%
\end{array}
\right\vert \\
D_{1}  &  =\left(  -1\right)  ^{\frac{3\cdot2}{2}}\left\vert
\begin{array}
[c]{ccccc}%
a_{4} & a_{3} & a_{2} & a_{1} & a_{0}\\
& a_{4} & a_{3} & a_{2} & a_{1}\\
4a_{4} & 3a_{3} & 2a_{2} & 1a_{1} & \\
& 4a_{4} & 3a_{3} & 2a_{2} & 1a_{1}\\
&  & 4a_{4} & 3a_{3} & 2a_{2}%
\end{array}
\right\vert \\
D_{2}  &  =\left(  -1\right)  ^{\frac{2\cdot1}{2}}\left\vert
\begin{array}
[c]{ccc}%
a_{4} & a_{3} & a_{2}\\
4a_{4} & 3a_{3} & 2a_{2}\\
& 4a_{4} & 3a_{3}%
\end{array}
\right\vert \\
D_{3}  &  =\left(  -1\right)  ^{\frac{1\cdot0}{2}}\left\vert
\begin{array}
[c]{c}%
4a_{4}%
\end{array}
\right\vert
\end{align*}

\end{example}

\begin{theorem}
\label{thm:sdis} For $k=1,\ldots,n-1$, we have
\[
Z_{k}=V\left( \{ D_{0},\ldots,D_{n-k-1}\}\right)  .
\]
Hence $Z_{n-1}\supsetneq Z_{n-2}\supsetneq\cdots\supsetneq Z_{1}$ is the
\emph{subdiscriminant}-based stratification of the discriminant
hypersurface~$Z_{n-1}$.
\end{theorem}

\section{Order-based Stratification}

\label{sec:order}

\subsection{Main result}

For the sake of simplicity, we use the following order-related notations in
this section.

\begin{notation}
For $P\in\mathbb{C}[u_{1},\ldots,u_{p}]$ and $\gamma\in\mathbb{C}^{p}$,

\begin{itemize}
\item $\partial_{(k_{1},\ldots,k_{p})}P:=\dfrac{\partial^{k}P}{\partial
u_{1}^{k_{1}}\cdots\partial u_{p}^{k_{p}}}$ where $k_{i}\in\mathbb{N}$ and
$k=k_{1}+\cdots+k_{p}$;

\item $\operatorname*{ord}P(\gamma):=$ the order of $P$ at $\gamma$.
\end{itemize}
\end{notation}

\begin{theorem}
\label{thm:order}For $k=1,\ldots,n-1$, we have $Z_{k}=O_{k}$ where%
\[
O_{k}=\left\{  \gamma\in\mathbb{C}^{n}:\operatorname*{ord}D(\gamma)\geq
n-k\right\}  .
\]
Equivalently,
\[
Z_{k}\ =\ V(\{\partial_{\delta}D\ :\ \left\vert \delta\right\vert
\ \leq\ n-k-1\}).
\]
Hence $Z_{n-1}\supsetneq Z_{n-2}\supsetneq\cdots\supsetneq Z_{1}$ is the
\emph{order}-based stratification of the discriminant hypersurface $Z_{n-1}$.
\end{theorem}

\begin{example}
\label{ex:order} Let $n=4$ and $k=2$. Then $F=x^{4}+a_{3}x^{3}+a_{2}%
x^{2}+a_{1}x+a_{0}$. For the polynomial $P_{2} $ in Example \ref{ex:quartic},
i.e., $\gamma=(4,0,4,0)$, it is obvious that $P_{2} $ has two distinct roots
and thus $\gamma\in Z_{2}$. By Theorem \ref{thm:order}, we have $\gamma\in
O_{2}$. In other words, $\operatorname*{ord}D(\gamma)\ge4-2=2.$

In fact, we can verify $\gamma\in O_{2}$ by calculation. Recall that
\begin{align*}
D =  &  \,256a_{0}^{3}-27a_{1}^{4}-128a_{0}^{2}a_{2}^{2}-192a_{0}^{2}%
a_{1}a_{3}+144a_{0}a_{1}^{2}a_{2} -4a_{1}^{2}a_{2}^{3}+18a_{1}^{3}a_{2}%
a_{3}-6a_{0}a_{1}^{2}a_{3}^{2}\\
&  +16a_{0}a_{2}^{4}+144a_{0}^{2}a_{2}a_{3}^{2}-80a_{0}a_{1}a_{2}^{2}%
a_{3}+18a_{0}a_{1}a_{2}a_{3}^{3} -4a_{1}^{3}a_{3}^{3}-4a_{0}a_{2}^{3}a_{3}%
^{2}-27a_{0}^{2}a_{3}^{4}+a_{1}^{2}a_{2}^{2}a_{3}^{2}.
\end{align*}
Taking the first partial derivative of $D$ with respect to $a_{i}$ for
$i=0,1,2,3$, we obtain
\begin{align*}
\dfrac{\partial D}{\partial a_{0}}=  &  -54\,a_{{0}}{a_{{3}}^{4}}+18\,a_{{1}%
}a_{{2}}{a_{{3}}^{3}}-4\,{a_{{2}}^{3}} {a_{{3}}^{2}}+288\,a_{{0}}a_{{2}%
}{a_{{3}}^{2}}-6\,{a_{{1}}^{2}}{a_{{3}}^{2}}-80\,a_{{1}}{a_{{2}}^{2}}a_{{3}%
}+16\,{a_{{2}}^{4}}-384\,a_{{0 }}a_{{1}}a_{{3}}\\
&  -256\,a_{{0}}{a_{{2}}^{2}}+144\,{a_{{1}}^{2}}a_{{2}}+ 768\,{a_{{0}}^{2}},\\
\dfrac{\partial D}{\partial a_{1}}=  &  \,18\,a_{{0}}a_{{2}}{a_{{3}}^{3}%
}-12\,{a_{{1}}^{2}}{a_{{3}}^{3}}+2\,a_{{\ 1}}{a_{{2}}^{2}}{a_{{3}}^{2}%
}-12\,a_{{0}}a_{{1}}{a_{{3}}^{2}}-80\,a_{{0 }}{a_{{2}}^{2}}a_{{3}}%
+54\,{a_{{1}}^{2}}a_{{2}}a_{{3}}-8\,a_{{1}}{a_{{\ 2}}^{3}}-192\,{a_{{0}}^{2}%
}a_{{3}}\\
&  +288\,a_{{0}}a_{{1}}a_{{2}}-108\,{a_{{1}}^{3}},\\
\dfrac{\partial D}{\partial a_{2}}=  &  \,18\,a_{{0}}a_{{1}}{a_{{3}}^{3}%
}-12\,a_{{0}}{a_{{2}}^{2}}{a_{{3}}^{2}}+ 2\,{a_{{1}}^{2}}a_{{2}}{a_{{3}}^{2}%
}+144\,{a_{{0}}^{2}}{a_{{3}}^{2}}- 160\,a_{{0}}a_{{1}}a_{{2}}a_{{3}%
}+64\,a_{{0}}{a_{{2}}^{3}}+18\,{a_{{1} }^{3}}a_{{3}}-12\,{a_{{1}}^{2}}{a_{{2}%
}^{2}}\\
&  -256\,{a_{{0}}^{2}}a_{{2}} +144\,a_{{0}}{a_{{1}}^{2}},\\
\dfrac{\partial D}{\partial a_{3}}=  &  \,18\,a_{{0}}a_{{1}}{a_{{3}}^{3}%
}-12\,a_{{0}}{a_{{2}}^{2}}{a_{{3}}^{2}}+ 2\,{a_{{1}}^{2}}a_{{2}}{a_{{3}}^{2}%
}+144\,{a_{{0}}^{2}}{a_{{3}}^{2}}- 160\,a_{{0}}a_{{1}}a_{{2}}a_{{3}%
}+64\,a_{{0}}{a_{{2}}^{3}}+18\,{a_{{1} }^{3}}a_{{3}}-12\,{a_{{1}}^{2}}{a_{{2}%
}^{2}}\\
&  -256\,{a_{{0}}^{2}}a_{{2}} +144\,a_{{0}}{a_{{1}}^{2}}.
\end{align*}
It can be verified that
\[
D(\gamma)=\dfrac{\partial D}{\partial a_{0}}(\gamma)=\dfrac{\partial
D}{\partial a_{1}}(\gamma)=\dfrac{\partial D}{\partial a_{2}}(\gamma
)=\dfrac{\partial D}{\partial a_{3}}(\gamma)=0
\]
and
\[
\dfrac{\partial^{2} D}{\partial a_{0}\partial a_{2}}(\gamma)=\left(
18\,a_{{1}}{a_{{3}}^{3}}-12\,{a_{{2}}^{2}}{a_{{3}}^{2}}+288\,a_{{0}}{a _{{3}%
}^{2}}-160\,a_{{1}}a_{{2}}a_{{3}}+64\,{a_{{2}}^{3}}-512\,a_{{0}}a _{{2}%
}+144\,{a_{{1}}^{2}} \right)  \big|_{a=\gamma}=-4096\ne0.
\]
Hence $\operatorname*{ord}D(\gamma)=2$, indicating that $\gamma\in O_{2} $.
\end{example}

\subsection{Proof of Theorem \ref{thm:order}}

\noindent Now we will prove Theorem \ref{thm:order}. For this we need the
following lemma.

\begin{lemma}
\label{lem:condition_partial_derivative}
The polynomial $F_{\gamma}$ has at most $m$ distinct roots if and only if
\begin{equation}
\dfrac{\partial^{0}D}{\partial a_{0}^{0}}(\gamma)\ =\ \dfrac{\partial^{1}%
D}{\partial a_{0}^{1}}(\gamma)\ \ =\ \ \cdots\ \ =\ \ \dfrac{\partial
^{n-m-1}D}{\partial a_{0}^{n-m-1}}(\gamma)\ \ =\ \ 0.
\label{eq:derivative_constant}%
\end{equation}

\end{lemma}

\begin{proof}
First, we derive an equivalent condition for Equation
\eqref{eq:derivative_constant}. By Taylor expansion,
\[
\operatorname{res}_{x}(F+t,F^{\prime})=D(a_{0}+t,a_{1},\ldots,a_{n-1}
)=\sum_{i\geq0}t^{i}\cdot\frac{1}{i!}\frac{\partial^{i}}{\partial a_{0}^{i}%
}D.
\]
Therefore, Equation \eqref{eq:derivative_constant} is equivalent to
$t^{n-m}\mid\operatorname{res}_{x}(F_{\gamma} +t,F_{\gamma}^{\prime})$.

Let $G_{\gamma}=\gcd\left(  F_{\gamma},F_{\gamma}^{\prime}\right)  $ and
$d=\deg G_{\gamma}$. Assume $F_{\gamma}^{\prime}=G_{\gamma}\overline
{F_{\gamma}^{\prime}}$. It is easy to see that
\begin{equation}
\operatorname*{res}\nolimits_{x}\left(  F_{\gamma}+t,F_{\gamma}^{\prime
}\right)  =\operatorname*{res}\nolimits_{x}\left(  F_{\gamma}+t,G_{\gamma
}\overline{F_{\gamma}^{\prime}}\right)  =\operatorname*{res}\nolimits_{x}
\left(  F_{\gamma}+t,G_{\gamma}\right)  \ \operatorname*{res}\nolimits_{x}
\left(  F_{\gamma}+t,\overline{F_{\gamma}^{\prime}}\right)  .
\label{eq:factor_decomp_res}%
\end{equation}

Let $\alpha_{1},\ldots,\alpha_{d}$ be the complex roots of $G_{\gamma}$.
Obviously, $F_{\gamma}(\alpha_{i})=0$. Hence
\begin{equation}
\operatorname*{res}\nolimits_{x}\left(  F_{\gamma}+t,G_{\gamma}\right)
=(-1)^{nd}\cdot1^{n}\cdot\prod_{i=1}^{d}(F_{\gamma}+t)(\alpha_{i})=(-1)^{nd}
\prod_{i=1}^{d}t=(-1)^{nd}t^{d}. \label{eq:factor_decomp_cofac}%
\end{equation}

Now we are going to show that
\[
F_{\gamma} \text{ has at most } m \text{ distinct roots} \Longleftrightarrow
t^{n-m}\mid\operatorname{res}_{x}(F_{\gamma} +t,F_{\gamma}^{\prime}).
\]

\begin{enumerate}
\item[$\Longrightarrow$:] Combining Equations \eqref{eq:factor_decomp_res} and
\eqref{eq:factor_decomp_cofac}, we have
\begin{equation}
\operatorname*{res}\nolimits_{x}\left(  F_{\gamma}+t,F_{\gamma}^{\prime
}\right)  =(-1)^{nd}t^{d}\ \operatorname*{res}\nolimits_{x}\left(  F_{\gamma
}+t,\overline{F_{\gamma}^{\prime}}\right)  . \label{eq:res_F_cofac(F')}%
\end{equation}

Therefore, if $F_{\gamma}$ has at most $m$ distinct roots, then $n-m\leq
d=\deg G_{\gamma}$. By Equation \eqref{eq:res_F_cofac(F')}, we immediately
have $t^{n-m}\mid\operatorname{res}_{x}(F_{\gamma}+t,F^{\prime}_{\gamma})$ .

\item[$\Longleftarrow$:] 

The combination of the assumption $t^{n-m} \mid\operatorname{res}%
_{x}(F_{\gamma}+t,F^{\prime}_{\gamma})$ and Equation
\eqref{eq:res_F_cofac(F')} yields
\[
t^{n-m}\mid t^{d}\ \operatorname*{res}\nolimits_{x}\left(  F_{\gamma
}+t,\overline{F_{\gamma}^{\prime}}\right)  .
\]
So in order to prove $d\geq n-m$, it suffices to prove that $t\nmid$
$\operatorname*{res}\nolimits_{x}\left(  F_{\gamma}+t,\overline{F^{\prime
}_{\gamma}}\right)  $, or equivalently, $\operatorname*{res}\nolimits_{x}%
\left(  F_{\gamma},\overline{F_{\gamma}^{\prime}}\right)  \neq0$.

But this is almost trivial, as
\[
\operatorname*{res}\nolimits_{x}\left(  F_{\gamma},\overline{F_{\gamma
}^{\prime}}\right)  \neq0\ \Longleftrightarrow\ \gcd(F_{\gamma},\overline
{{F_{\gamma}^{\prime}}})=1 \ \Longleftrightarrow\ \gcd(G_{\gamma}%
,\overline{{F_{\gamma}^{\prime}}})=1,
\]
and $G_{\gamma},\overline{{F_{\gamma}^{\prime}}}$ are clearly coprime by
counting the multiplicity of common roots of $F_{\gamma}$ and $F_{\gamma
}^{\prime}$.
\end{enumerate}
\end{proof}

\begin{remark}
\label{rem:pdiff-a0} In Lemma \ref{lem:condition_partial_derivative}, we
obtain the following nice result:
\[
Z_{m}\ =\ V\left(\left\{  D,\dfrac{\partial D}{\partial a_{0}},\ldots,\dfrac
{\partial^{n-m-1}D}{\partial a_{0}^{n-m-1}}\right\}\right)  .
\]
This characterizes $Z_{m}$ in terms of the partial derivatives of $D$ with
respect to $a_{0}$. A curious reader might wonder what if $a_{1},a_{2}$ or
other coefficients were used instead of $a_{0}$.

In fact, similar relations hold, but slight adjustments have to be made in
both the statement and the proof. We prepare a satisfactory answer in the
Appendix for those with the greatest curiosity.

\end{remark}

Now we are ready to prove Theorem \ref{thm:order}.

\begin{proof}
[Proof of Theorem \ref{thm:order}] The proof will be divided into
several steps for easier reading.

\begin{enumerate}
\item It is easy to see that%
\[
\text{the order of }D\ \text{at }\gamma\ \text{with respect to}\ a\geq
n-k\ \ \Longleftrightarrow\ \ \gamma\in V(\{\partial_{\delta}\left(  D\right)
:\,0\leq|\delta|\leq n-k-1\})
\]
where $|\delta|$ is the shorthand notation of the sum of elements in $\delta$
hereinafter. Equivalently,
\[
O_{k}=V(\{\partial_{\delta}\left(  D\right)  |\,0\leq|\delta|\leq n-k-1\}).
\]
Next we will prove $Z_{k}=O_{k}$ by showing $Z_{k}\subseteq O_{k}$ and
$Z_{k}\supseteq O_{k}$.

\item $Z_{k}\subseteq O_{k}$.

\begin{enumerate}
\item Let $\gamma\in Z_{k}$. Then $F_{\gamma}$ has at most $k$ distinct roots.
Hence
\[
\deg\gcd(F_{\gamma},F_{\gamma}^{\prime})\geq n-k.
\]

\item Let $M_{\gamma}$ denote the Sylvester
matrix of $F_{\gamma}$ and $F_{\gamma}^{\prime}$ with respect to $x $. By
\cite[Chapter 1, Proposition 1]{2006_Apery_Jouanolou},
\[
\operatorname*{rank}M_{\gamma}=2n-1-\deg\gcd(F_{\gamma},F^{\prime}_{\gamma})\leq2n-1-(n-k)=n+k-1.
\]

\item Assume $\delta=(\delta_{0},\ldots,\delta_{n-1})$ satisfies $|\delta|\leq
n-k-1$. Let $e_{i}$ denote the $i$-th unit vector of length $n$. One can write
$\delta$ into the sum of $|\delta|$ unit vectors, i.e., $\delta=\sum
_{j=1}^{|\delta|}e_{s_{j}}$. Thus
\[
\partial_{\delta}D=\partial_{e_{s_{|\delta|}}}\left(  \cdots\partial
_{e_{s_{2}}}\left(  \partial_{e_{s_{1}}}D\right)  \right)  .
\]
Let $M$ be the Sylvester
matrix of $F$ and $F^{\prime}$ with respect to $x$. By \cite[Chapter 11, Section 1, Remark on Theorem 1]{2018_Adams_Essex},
\[
\partial_{e_{s_{1}}}D=\sum_{i=1}^{n}\det[M_{1},\ldots,M_{i-1},\hat{M}%
_{i},M_{i+1},\ldots,M_{n}],
\]
where $M_{j}$ is the $j$-th column of $M$, and $\hat{M}%
_{i}=\dfrac{\partial M_{i}}{\partial a_{{s_{1}}-1}}$ with the partial
derivative taken entrywise. Note that the entries of $M$ are linear with respect to $a_{i}$'s. Thus each $\hat{M}_{i}
$ is a constant vector. Taking the Laplace expansion along the column $\hat
{M}_{i}$ for each of the above summands, we can write $\partial_{e_{s_{1}}}D$
as a sum of minors of~$M$ with size $2n-2$.

Now apply the $\partial_{e_{s_{2}}}$ operator on each minor. Noting that the
entries of each minor are again linear with respect to $a_{i}$'s, we conclude
that $\partial_{e_{s_{2}}}\left(  \partial_{e_{s_{1}}}D\right)  $ is a sum of
minors of $M$ with size~$2n-3$.

Repeating the process, we can eventually write $\partial_{\delta}D$ as a sum
of minors of $M$ with size~$2n-1-|\delta|$.

\item From (b), we deduce that
\[
2n-1-|\delta|\geq2n-1-(n-k-1)=n+k>n+k-1\geq\operatorname*{rank}M_{\gamma}.
\]
Thus the partial derivative $(\partial_{\delta}D)(\gamma)$ is a sum of minors of
$M_{\gamma}$ with size greater than the
rank. Hence each of the summands is zero, which implies that $(\partial
_{\delta}D)(\gamma)=0$. Hence $\gamma\in O_{k}$.
\end{enumerate}

\item $Z_{k}\supseteq O_{k}$.

From Lemma \ref{lem:condition_partial_derivative}, we have
\[
Z_{k}\ =\ V\left( \left\{ D,\dfrac{\partial D}{\partial a_{0}},\ldots,\dfrac
{\partial^{n-k-1}D}{\partial a_{0}^{n-k-1}}\right\}\right)  .
\]
It immediately follows that $Z_{k}\supseteq O_{k}$.
\end{enumerate}
\end{proof}

\section{Smoothness-based Stratification}

\label{sec:smooth}

\begin{theorem}
\label{thm:singular_stratification} For $k=n-1,n-2,\ldots,2$, we have
\[
\operatorname*{Sing}Z_{k}=Z_{k-1}.
\]
Equivalently, for $k=1,\ldots,n-1,$ we have%
\[
Z_{k}=\ \underset{n-1-k\text{ iterations}}{\underbrace{\operatorname*{Sing}%
\cdots\operatorname*{Sing}}}\ V(D).
\]
Hence $Z_{n-1}\supsetneq Z_{n-2}\supsetneq\cdots\supsetneq Z_{1}$ is the
\emph{smoothness}-based stratification of the discriminant
hypersurface~$Z_{n-1}$.
\end{theorem}

Now we will prove Theorem \ref{thm:singular_stratification}. For this, we need
to investigate the irreducible decomposition of $Z_{k}$. Since a monic
polynomial of degree $n$ is determined by its non-leading coefficients, the
affine space $\mathbb{C}^{n}$ parameterizes the set of monic polynomials of
degree $n$.

\begin{definition}
Let $n$ be a positive integer and let $\mu=(\mu_{1},\ldots,\mu_{m})$ be a
partition of $n$.
We say a monic polynomial $P$ of degree $n$ lies in the \textbf{coincidence
root locus} $C_{\mu}$, if $P(x)$ splits into $m$ (not necessarily distinct)
linear factors raised to powers $\mu_{1},\ldots,\mu_{m}$ over $\mathbb{C}$:
\[
P(x)=\prod_{i=1}^{m}(x-x_{i})^{\mu_{i}}.
\]
In other words,
\[
C_{\mu}=\left\{  P\in\mathbb{C}^{n}\middle|P(x)=\prod_{i=1}^{m}(x-x_{i}%
)^{\mu_{i}}\right\}  .
\]

\end{definition}

\begin{example}
For example, let $P=x^{4}-2x^{2}+1$. Then
\[
P=(x-1)^{2}(x+1)^{2}=(x-1)(x-1)(x+1)^{2}=(x-1)(x-1)(x+1)(x+1).
\]
Therefore $P$ lies in $C_{(1,1,1,1)},C_{(1,1,2)}$ and $C_{(2,2)}$.

\end{example}

\begin{proposition}
[Folklore]\label{prop:irr-dec-of-subdisc-var} The set $C_{\mu}$ is an
irreducible closed subset of $\mathbb{C}^{n}$, of dimension $\#\mu$.
\end{proposition}

\begin{proof}
Consider the Vi\`{e}ta map $\nu$ sending roots $(x_{1},\ldots,x_{n})$ to the
coefficients $(a_{0},\ldots,a_{n-1})$ defined by the elementary symmetric
polynomials. The map $\nu$ is a finite flat surjective morphism from
$\mathbb{C}^{n}$ to $\mathbb{C}^{n}$, because the polynomial ring
$\mathbb{C}[x_{1},\ldots,x_{n}]$ is a free module of rank $n!$ over the
subring of symmetric polynomials \cite[Chapter IV, Section 6, No. 1, Theorem
1(c)]{2013_Bourbaki}. Then $C_{\mu}$ is the image of a linear subspace of
dimension $\#\mu$ under $\nu$, hence irreducible, closed and of dimension
$\#\mu$.
\end{proof}

\begin{example}
For example, the set $C_{(2,2)}\subseteq\mathbb{C}^{4}$ is the image of
$H=\left\{  (x_{1},\ldots,x_{4})\in\mathbb{C}^{4}\middle|x_{1}=x_{2}%
,x_{3}=x_{4}\right\}  $ under the Vi\`{e}ta map $\nu$.
Clearly $H\simeq\mathbb{C}^{2}$ is a two-dimensional irreducible closed set,
so is $C_{(2,2)}$.
\end{example}

\begin{corollary}
\label{cor:irr-dec-of-v-k} We have the following irreducible decomposition of
$Z_{k}$:
\[
Z_{k}=\bigcup_{\#\mu=k}C_{\mu}.
\]
Therefore, $Z_{k}$ is equidimensional and each $C_{\mu}$ is an irreducible
component of $Z_{k}$.
\end{corollary}

\begin{proof}
This immediately follows from Proposition \ref{prop:irr-dec-of-subdisc-var}
and the fact that a polynomial $P$ has at most $k$ distinct roots if and only
if $P$ lies in a coincident root locus $C_{\mu}$ for some $\#\mu=k$. Also
notice that $\dim C_{\mu}=\#\mu$.\bigskip
\end{proof}

\noindent Now we are ready to prove Theorem \ref{thm:singular_stratification}.

\begin{proof}
[Proof of Theorem \ref{thm:singular_stratification}] Let $\bigcup
_{\#\mu=k}C_{\mu}$ be the irreducible decomposition of the variety $Z_{k}$.
Then we have
\begin{equation}
\operatorname*{Sing}Z_{k}=\bigcup_{\substack{\mu\neq\mu^{\prime} \\\#\mu
=\#\mu^{\prime}=k}}(C_{\mu}\cap C_{\mu^{\prime}})\cup\bigcup_{\#\mu
=k}\operatorname*{Sing}C_{\mu} \label{eqn:singular locus}%
\end{equation}
by \cite[Exercise 9.6.11(c)]{2015_Cox_Little_OShea}.

\begin{itemize}
\item $\operatorname*{Sing}Z_{k}\subseteq Z_{k-1}$.\newline

It suffices to show that, for $\mu\neq\mu^{\prime}$ and $\#\mu=\#\mu^{\prime
}=k$, we have $C_{\mu}\cap C_{\mu^{\prime}}$ $\subseteq Z_{k-1}$ and
$\operatorname*{Sing}C_{\mu}$ $\subseteq Z_{k-1}$.

\begin{enumerate}
\item For every monic polynomial $P\in C_{\mu}\cap C_{\mu^{\prime}}$ with
$\#\mu=\#\mu^{\prime}=k$, we see that the (monic) linear factors of $P$ can be
grouped in two different ways, namely $\mu$ and $\mu^{\prime}$. Thus the
number of distinct roots of $P$ must be less than $k$. Therefore,
\[
C_{\mu}\cap C_{\mu^{\prime}}\subseteq Z_{k-1}.
\]

\item In {\cite[Theorem 5.4]{2003_Chipalkatti}}, Chipalkatti showed that for
every $\mu$ the set $\operatorname*{Sing}C_{\mu}$ is a union of certain
$C_{\lambda}$'s where each $\lambda$ is a coarsening of $\mu$. Note that
$\#\lambda<k$ and $C_{\lambda}\subseteq Z_{\#\lambda}\subseteq Z_{k-1}$.
Hence
\[
\operatorname*{Sing}C_{\mu}\subseteq Z_{k-1}.
\]

\end{enumerate}

\item $\operatorname*{Sing}Z_{k}\supseteq Z_{k-1}$. \newline Now we will show
that all irreducible components of $Z_{k-1}$ are contained in
$\operatorname*{Sing}Z_{k}$. This is done by analyzing the irreducible
components of $Z_{k-1}=\bigcup_{\#\lambda=k-1}C_{\lambda}$ by the largest
number in $\lambda=(\lambda_{1},\ldots,\lambda_{k-1})$.

\begin{enumerate}
\item Suppose $\max_{i}\lambda_{i}\geq4$. Without loss of generality, we may
assume that $\lambda_{1}\geq4$. Then $\lambda$ can be refined to two different
partitions with length $k$, e.g., $\mu=(1,\lambda_{1}-1,\lambda_{2}%
,\ldots,\lambda_{k-1})$ and $\mu^{\prime}=(2,\lambda_{1}-2,\lambda_{2}%
,\ldots,\lambda_{k-1})$. So $C_{\lambda}\subseteq(C_{\mu}\cap C_{\mu^{\prime}%
})\subseteq\operatorname*{Sing}Z_{k}$.

\item Suppose $\max_{i}\lambda_{i}=3$. Let
\[
\lambda=(\underbrace{3,\ldots,3}_{\ell},\underbrace{2,\ldots,2}_{\ell^{\prime
}},\underbrace{1,\ldots,1}_{\ell^{\prime\prime}})
\]
and
\[
\mu=(\underbrace{3,\ldots,3}_{\ell-1},\underbrace{2,\ldots,2}_{\ell^{\prime
}+1},\underbrace{1,\ldots,1}_{\ell^{\prime\prime}+1}).
\]
Note that $\lambda$ is a coarsening of $\mu$ by merging $2$ and $1$ from $\mu$
into $3$. From {\cite[Theorem 5.4]{2003_Chipalkatti}}, we have $C_{\lambda
}\subseteq\operatorname*{Sing}C_{\mu}$. Recalling (\ref{eqn:singular locus}),
we see that $\operatorname*{Sing}C_{\mu}\subseteq\operatorname*{Sing}Z_{k}$.
Hence we have
\[
C_{\lambda}\subseteq\operatorname*{Sing}Z_{k}.
\]

\item Suppose $\max_{i}\lambda_{i}=2$. Let
\[
\lambda=(\underbrace{2,\ldots,2}_{\ell},\underbrace{1,\ldots,1}_{\ell^{\prime
}})
\]
and
\[
\mu=(\underbrace{2,\ldots,2}_{\ell-1},\underbrace{1,\ldots,1}_{\ell^{\prime
}+2}).
\]
Note that $\lambda$ is a coarsening of $\mu$ by merging two copies of $1$ from
$\mu$ into $2$, where $2$ is already in $\mu$. From {\cite[Theorem
5.4]{2003_Chipalkatti}}, we have $C_{\lambda}\subseteq\operatorname*{Sing}%
C_{\mu}$. Recalling (\ref{eqn:singular locus}), we see that
$\operatorname*{Sing}C_{\mu}\subseteq\operatorname*{Sing}Z_{k}$ Hence we have
\[
C_{\lambda}\subseteq\operatorname*{Sing}Z_{k}.
\]

\end{enumerate}

Therefore $\operatorname*{Sing}Z_{k}\supseteq Z_{k-1}$.
\end{itemize}

So we have the claimed equality $\operatorname*{Sing}Z_{k}=Z_{k-1}$. Notice
that $Z_{1}$ is isomorphic to the affine line because it is the image of the
map
\[
t\mapsto\left(  (-1)^{1}\binom{n}{1}t,(-1)^{2}\binom{n}{2}t^{2},\ldots
,(-1)^{n}\binom{n}{n}t^{n}\right)
\]
(recall Vieta's formulas). Thus $Z_{1}$ is smooth and $Z_{n-1}\supsetneq
Z_{n-2}\supsetneq\cdots\supsetneq Z_{1}$ is indeed the smooth stratification
of $Z_{n-1}$.
\end{proof}

\section{Conclusion}

\label{sec:conclusion}

\paragraph{Summary:}

In this paper, we considered stratification of the discriminant hypersurface:
\[
Z_{n-1}\supsetneq Z_{n-2}\supsetneq\cdots\supsetneq Z_{1}%
\]
where $Z_{k}$ is the set of all monic univariate polynomials of degree $n$
with at most $k$ distinct roots. We reviewed a classical approach to describe
$Z_{k}$ and then proposed two new approaches.

\begin{enumerate}
\item In Section~\ref{sec:sdis} (Theorem \ref{thm:sdis}), we briefly reviewed
a classical description of $Z_{k}$ via subdiscriminants $D_{i}$'s. Precisely
put,
\[
Z_{k}\ \ =\ \ V\left( \{D_{0},D_{1},\ldots,D_{n-k-1}\}\right)  .\ \ \ \ \ \
\]

\item In Section~\ref{sec:order} (Theorem \ref{thm:order}), we gave a new
description of $Z_{k}$ via partial derivatives of the discriminant~$D$ (that
is $D_{0}$). Precisely put,%
\[
Z_{k} \ = \ V( \{\partial_{\delta}D \ : \ \left\vert \delta\right\vert
\ \leq\ n-k-1\}).
\]

\item In Section~\ref{sec:smooth} (Theorem \ref{thm:singular_stratification}),
we gave another new description of $Z_{k}$ via repeated singular locus of the
discriminant $D$. Precisely put%
\[
Z_{k}=\underset{n-1-k\text{ iterations}}{\underbrace{\operatorname*{Sing}%
\cdots\operatorname*{Sing}}}\ V\left(  D\right)  .\ \ \ \ \ \ \ \ \ \ \ \ \ \
\]

\end{enumerate}

\paragraph{Future works:}

\begin{itemize}
\item Note that the first two descriptions of $Z_{k}$ imply that%
\[
V(D_{0},D_{1},\ldots,D_{n-k-1})\ =\ V(\langle\partial_{\delta}D\ \colon
\ |\delta|\ \leq\ n-k-1\rangle),
\]
or equivalently,%
\[
\sqrt{\langle D_{0},D_{1},\ldots,D_{n-k-1}\rangle}\ =\ \sqrt{\langle
\partial_{\delta}D\ :\ |\delta|\ \leq\ n-k-1\ \rangle}%
\]
However, numerical evidences suggest that
\begin{align*}
&  \langle D_{0},D_{1},\ldots,D_{n-k-1}\rangle\ \neq\ \langle\partial_{\delta
}D\ :\ |\delta|\ \leq\ n-k-1\ \rangle\\
&  \langle D_{0},D_{1},\ldots,D_{n-k-1}\rangle\ \ \ \ \ \ \ \ \text{is not
radical}\\
\  &  \langle\partial_{\delta}D\ :\ |\delta|\ \leq\ n-k-1\ \rangle\ \text{is
not radical}%
\end{align*}
Thus a question naturally arises:

\hspace{10em}\emph{Is there a \textquotedblleft good\textquotedblright\ set of
generators for their radicals?}

Knowing those would allow us to apply the Jacobian criterion to write down the
defining polynomials of the singular locus of $Z_{k}$ as we did in the third
approach.

\item Recall the relation from the above%
\[
V(\{ D_{0},D_{1},\ldots,D_{n-k-1}\} )\ =\ V(\{\partial_{\delta}D\ \colon
\ |\delta|\ \leq\ n-k-1\})
\]
Hilbert's Nullstellensatz implies that
\[
\underset{s}{\exists}D_{k}^{s}\in\langle\partial_{\delta}D\ :\ |\delta
|\ \leq\ n-k-1\ \rangle
\]
Another question naturally arises: \emph{What is the smallest number for $s$?}

\smallskip Based on some computations, we conjecture that the smallest
exponent is $k+1$. 
The claim has been verified  for small $n$ and $k$ using a Maple program that can be downloaded from the following 
site: 
\begin{center}
\vskip-1em
\url{https://github.com/JYangMATH/DiscriminantHypersurface}
\end{center}
\vskip-1em
However, it is still unclear how to prove/disprove the claim for the
general case.

\item In this paper, we partition the discriminant hypersurface by the number
of distinct roots. Recall that the number of distinct roots can be captured
via subdiscriminants and their variants
\cite{2009_Sturm,1996_Yang_Hou_Zeng,1998_Gonzalez_Vega_Recio_Lombardi_Roy} and
(signed) subdiscriminants can also be used to retrieve the information of real
roots. Then a natural question arises: can we partition the discriminant
hypersurface by the number of distinct roots in the real space? This question
is worthy of further investigation in the future.

\item
Finally the most  significant generalization of the current work is to extend the theory of the discriminant hypersurface from the single univariate polynomial case to the  multivariate polynomial system case.
For this,  one would need to revisit the works on multivariate discriminants
\cite{
1994_Gelfand_Kapranov_Zelevinsky,
2013_Esterov,
2016_Han,
2023_Dickenstein_DiRocco_Morrison}
and/or on multivariate subresultants~\cite{1995_Chardin}
in view of stratifying  the discriminant hypersurface in terms of number of distinct complex/real solutions.
\end{itemize}

\section*{Appendix}
The main result of the paper does \emph{not} depend on the result hereinafter.
This appendix is for the readers who want to learn more about details behind Remark \ref{rem:pdiff-a0}. 

In Lemma \ref{lem:condition_partial_derivative}, we show that $Z_k$ can be characterized as the vanishing set of $D$ and first $n-k-1$ partial derivatives of $D$ with respect to $a_0$. This is no longer true if we replace $a_0$ with any other coefficients. However it is not very far from being true if we are willing to make a small adjustment.
\begin{proposition}\label{prop:Zk_unmixed_partial_derivative}
For $ 1\leq s\leq n-1$ and $1\le k\le n-1$, we have 
$$Z_k=\overline{ V\left(\left\{D,\dfrac{\partial D}{\partial a_{s}},\ldots,\dfrac{\partial^{n-k-1}D}{\partial a_{s}^{n-k-1}}\right\}\right)\big\backslash V\left(  a_{0}\right)}.$$
\end{proposition}

In order to prove Proposition \ref{prop:Zk_unmixed_partial_derivative}, we need the following lemma.

\begin{lemma}\label{lem:unmixed_derivative}
        Suppose $P=x^n+c_{n-1} x^{n-1}+\cdots+c_0$ is a monic polynomial that splits as $P=\prod_{i=1}^m(x-r_i)^{\mu_i}$, where $r_1,\ldots,r_m$ are distinct complex roots. Let $1\leq s\leq n-1$ be a positive integer. Then
        \begin{enumerate}
                \item $D(c_0,\ldots,c_{n-1})=\dfrac{\partial D}{\partial a_s}(c_0,\ldots,c_{n-1})=\cdots=\dfrac{\partial^{n-m-1} D}{\partial a_s^{n-m-1}}(c_0,\ldots,c_{n-1})=0$;
                \item $\dfrac{\partial^{n-m} D}{\partial a_s^{n-m}}(c_0,\ldots,c_{n-1})=0$ if and only if $0$ is a multiple root of $P$.
        \end{enumerate}
\end{lemma}
\begin{proof}
        Recall that $F=x^n+a_{n-1} x^{n-1}+\cdots+a_0$ is a general monic univariate polynomial of degree $n$ and $D$ is the discriminant of $F$. Now we introduce a fresh variable $t$ and consider $F_s=F+t x^s$. The discriminant of $F_s$ will simply be $D_s=D(a_0,\ldots,a_{s-1},a_s+t,a_{s+1},\ldots,a_{n-1})\in \mathbb{C}[a_0,\ldots,a_{n-1},t]$. By Taylor's formula, we have 
        $$D_s=\sum_{i=0}^{+\infty} \left(\frac{1}{i!}\dfrac{\partial^i D}{\partial^i_{a_s}}\right)t^i.$$
        Notice that this is actually a finite sum.
        
        To show the lemma, it suffices to show that 
        \begin{enumerate}
                \item $t^{n-m}\big|D_s(c_0,\ldots,c_{n-1};t)$, and
                \item $t^{n-m+1}\big|D_s(c_0,\ldots,c_{n-1};t)$ if and only if $0$ is a multiple root of $P$. %
        \end{enumerate} 
        
        Write $Q=\prod_{i=1}^m (x-r_i)^{\mu_i-1}=\gcd(P,P^{\prime})$.
        Notice that 
        \begin{align*}
                D_s(c_0,\ldots,c_{n-1};t)&=\operatorname*{res}\nolimits_{x}\left(P+t x^s,P^{\prime}+stx^{s-1}\right)\\
                                                                 &=\operatorname*{res}\nolimits_{x}\left(P+t x^s-\frac{1}{s}x\left(P^{\prime}+stx^{s-1}\right),P^{\prime}+stx^{s-1}\right)\\
                                                                 &=\operatorname*{res}\nolimits_{x}\left(P-\frac{1}{s}xP^{\prime},P^{\prime}+stx^{s-1}\right)\\
                                                                 &=\operatorname*{res}\nolimits_{x}\left(Q\frac{P-\frac{1}{s}xP^{\prime}}{Q},P^{\prime}+stx^{s-1}\right)\\
                                                                 &=\operatorname*{res}\nolimits_{x}\left(Q,P^{\prime}+stx^{s-1}\right)\operatorname*{res}\nolimits_{x}\left(\frac{P-\frac{1}{s}xP^{\prime}}{Q},P^{\prime}+stx^{s-1}\right)\\
                                                                 &=\operatorname*{res}\nolimits_{x}\left(Q,stx^{s-1}\right) \operatorname*{res}\nolimits_{x}\left(\frac{P-\frac{1}{s}xP^{\prime}}{Q},P^{\prime}+stx^{s-1}\right).\\                                                          
        \end{align*}
        Now we claim the following for the first and the second factor respectively:
        \begin{enumerate}
                \item $t^{n-m}\big|\operatorname*{res}\nolimits_{x}\left(Q,stx^{s-1}\right)$, and
                \item $t\big|\operatorname*{res}\nolimits_{x}\left(\frac{P-\frac{1}{s}xP^{\prime}}{Q},P^{\prime}+stx^{s-1}\right)$ if and only if $0$ is a multiple root of $P$.
        \end{enumerate}
                
        For the first claim, we have 
        \begin{align*}
                \operatorname*{res}\nolimits_{x}\left(Q,stx^{s-1}\right) & = \operatorname*{res}\nolimits_{x}\left(Q,t\right) \operatorname*{res}\nolimits_{x}\left(Q,sx^{s-1}\right)\\
                & = \operatorname*{res}\nolimits_{x}\left(\prod_{1\leq i\leq m} (x-r_i)^{\mu_i-1},t\right) \operatorname*{res}\nolimits_{x}\left(\prod_{1\leq i\leq m} (x-r_i)^{\mu_i-1},sx^{s-1}\right)\\
                & = t^{n-m}\cdot s^{n-m}\cdot \prod_{1\leq i \leq m}(r_i-0)^{(\mu_i-1)(s-1)} \\
                & = t^{n-m}\cdot s^{n-m}\cdot \prod_{1\leq i \leq m}r_i^{(\mu_i-1)(s-1)}.\\
        \end{align*}
        Therefore the first factor is indeed a multiple of $t^{n-m}$. The first part of our lemma follows from this claim. Notice also that if $0$ is not a multiple root of $P$, then $\operatorname*{res}\nolimits_{x}\left(Q,stx^{s-1}\right)=\alpha t^{n-m}$ for some non-zero constant $\alpha$.
        
        For the second claim, we will compute the constant term of $\operatorname*{res}\nolimits_{x}\left(\frac{P-\frac{1}{s}xP^{\prime}}{Q},P^{\prime}+stx^{s-1}\right)$, which is 
        \begin{align*}
                \left.\operatorname*{res}\nolimits_{x}\left(\frac{P-\frac{1}{s}xP^{\prime}}{Q},P^{\prime}+stx^{s-1}\right)\right|_{t=0}&= \operatorname*{res}\nolimits_{x}\left(\frac{P-\frac{1}{s}xP^{\prime}}{Q},P^{\prime}\right) \\
                &=\operatorname*{res}\nolimits_{x}\left(\frac{P-\frac{1}{s}xP^{\prime}}{Q},Q\frac{P^{\prime}}{Q}\right)\\
                &=\operatorname*{res}\nolimits_{x}\left(\frac{P-\frac{1}{s}xP^{\prime}}{Q},Q\right)\operatorname*{res}\nolimits_{x}\left(\frac{P-\frac{1}{s}xP^{\prime}}{Q},\frac{P^{\prime}}{Q}\right)\\
                &=\operatorname*{res}\nolimits_{x}\left(\frac{P-\frac{1}{s}xP^{\prime}}{Q},Q\right)\operatorname*{res}\nolimits_{x}\left(\frac{P}{Q},\frac{P^{\prime}}{Q}\right).\\
        \end{align*}
        Notice that $\operatorname*{res}\nolimits_{x}\left(\frac{P}{Q},\frac{P^{\prime}}{Q}\right)$ is clearly non-zero, as these polynomials are coprime. So it suffices to show that $\operatorname*{res}\nolimits_{x}\left(\frac{P-\frac{1}{s}xP^{\prime}}{Q},Q\right)$ vanishes if and only if $0$ is a multiple root of $P$. A direct computation shows
        \begin{align*}
                \operatorname*{res}\nolimits_{x}\left(\frac{P-\frac{1}{s}xP^{\prime}}{Q},Q\right) &= \operatorname*{res}\nolimits_{x}\left(\frac{\prod\limits_{i=1}^m(x-r_i)^{\mu_i}-\frac{1}{s}x\sum\limits_{i=1}^m\mu_i(x-r_i)^{\mu_i-1}\prod\limits_{j\ne
                                 i}(x-r_j)^{\mu_j}}{\prod\limits_{i=1}^m(x-r_i)^{\mu_i-1}},\prod_{i=1}^m(x-r_i)^{\mu_i-1}\right)\\
                &= \operatorname*{res}\nolimits_{x}\left(\prod\limits_{i=1}^m(x-r_i)-\frac{1}{s}x\sum\limits_{i=1}^m\mu_i\prod\limits_{j\ne i}(x-r_j),\prod_{i=1}^m(x-r_i)^{\mu_i-1}\right)\\
                &= (-1)^{m(n-m)}\prod_{k=1}^m \left(\prod\limits_{i=1}^m(r_k-r_i)-\frac{1}{s}r_k\sum\limits_{i=1}^m\mu_i\prod\limits_{j\ne i}(r_k-r_j)\right)^{\mu_k-1}\\
                &= (-1)^{m(n-m)}\prod_{k=1}^m \left(-\frac{1}{s}r_k\mu_k\prod\limits_{j\ne k}(r_k-r_j)\right)^{\mu_k-1}.\\
        \end{align*}
        The last line is obtained by observing that any product inside the biggest parenthesis in the previous line would be zero if $i=k$ or $j=k$.
        
        It is now obvious that $\operatorname*{res}\nolimits_{x}\left(\frac{P-\frac{1}{s}xP^{\prime}}{Q},Q\right)$ is zero if and only if $r_k=0$ and $\mu_k>1$, i.e.\ $0$ is a multiple root of $P$. This shows our second claim and the second part of the lemma follows.
\end{proof}

Now we are ready to prove Proposition \ref{prop:Zk_unmixed_partial_derivative}.

\begin{proof}[Proof of Proposition \ref{prop:Zk_unmixed_partial_derivative}]
We divide the proof into the following four steps:
\begin{enumerate}
\item 
Show $Z_{k}\subseteq V\left(\left\{D,\dfrac{\partial D}{\partial a_{s}},\ldots,\dfrac{\partial^{n-k-1}D}{\partial a_{s}^{n-k-1}}\right\}\right)$.

Suppose $\gamma\in Z_{k}$ and the number of distinct roots of $F_{\gamma}$ is $m$. Then $m\le k$. By Lemma \ref{lem:unmixed_derivative},  we have
\[(\partial_{a_{s}}^{\ell}D)(\gamma)=0\ \ \text{for}\ s\ge1, \ell=0,\ldots,n-m-1.\]
Since $m\le k$, we have $n-k-1\le n-m-1$ and thus 
\[(\partial_{a_{s}}^{0}D)(\gamma)\ =\ \cdots\ =\ (\partial_{a_{s}}^{n-k-1}D)(\gamma)\ =\ 0.\]
In other words, $\gamma\in V\left(\left\{D,\dfrac{\partial D}{\partial a_{s}},\ldots,\dfrac{\partial^{n-k-1}D}{\partial a_{s}^{n-k-1}}\right\}\right)$.

\item 
Show $ Z_{k}\supseteq V\left(\left\{D,\dfrac{\partial D}{\partial a_{s}},\ldots,\dfrac{\partial^{n-k-1}D}{\partial a_{s}^{n-k-1}}\right\}\right)\big\backslash V(a_0)$.

We show by contradiction. Let $\gamma \in V\left(\left\{D,\dfrac{\partial D}{\partial a_{s}},\ldots,\dfrac{\partial^{n-k-1}D}{\partial a_{s}^{n-k-1}}\right\}\right)\big\backslash V(a_0)$ but suppose that $\gamma \notin Z_k$. Let $m>k$ be the number of distinct roots of $F_{\gamma}$. By Lemma \ref{lem:unmixed_derivative}, $\dfrac{\partial^{n-m}D}{\partial a_s}$ does not vanish at $\gamma$ since $\gamma_0\ne 0$ and $0$ is not even a root of $F_\gamma$. However, this contradicts with our choice of $\gamma$. So we must have $\gamma \in Z_k$.

\item Show $ Z_{k}\backslash V(a_0)=V\left(\left\{D,\dfrac{\partial D}{\partial a_{s}},\ldots,\dfrac{\partial^{n-k-1}D}{\partial a_{s}^{n-k-1}}\right\} \right)\big\backslash V(a_0)$.

Combining 1 and 2, we have
\[V\left(\left\{D,\dfrac{\partial D}{\partial a_{s}},\ldots,\dfrac{\partial^{n-k-1}D}{\partial a_{s}^{n-k-1}}\right\}\right)\big\backslash V(a_0)\ \subseteq\ Z_{k}\ \subseteq\ V\left(\left\{D,\dfrac{\partial D}{\partial a_{s}},\ldots,\dfrac{\partial^{n-k-1}D}{\partial a_{s}^{n-k-1}}\right\}\right).\]
Removing $V(a_0)$ from the above three sets, we obtain that the desired relationship.

\item Show $Z_{k}=\overline{Z_k \backslash V(a_0)} $.

Note that none of the irreducible components of $Z_k$ is contained in $V(a_0)$. So if $W$ is an irreducible component of $Z_k$, $\overline{W\backslash V(a_0)}=W$. It follows that
        \begin{align*}
                \overline{Z_k \backslash V(a_0)} &= \overline{ \left(\bigcup_{{W\in \mathop{\mathrm{Comp}}(Z_k) }} W \right)\backslash V(a_0)}\\
                &=\overline{ \bigcup_{W\in \mathop{\mathrm{Comp}}(Z_k) }( W \backslash V(a_0))} \\
                &= \bigcup_{W\in \mathop{\mathrm{Comp}}(Z_k) } \overline{ W \backslash V(a_0)} \\
                &= \bigcup_{W\in \mathop{\mathrm{Comp}}(Z_k) } W \\
                &= Z_k.&
        \end{align*}
\end{enumerate}
Finally, the proof is completed.
\end{proof}

\bigskip\noindent\textbf{Acknowledgements.} 
Rizeng Chen's work was supported by National Key R\&D Program of China (Nos.: 2022YFA1005102 and 2024YFA1014003).
Hoon Hong's work was supported by
National Science Foundations of USA (Grant Nos: CCF 2212461 and CCF 2331401).
Jing Yang's work was supported by National Natural Science Foundation of China
(Grant Nos.: 12261010, 12326363 and  12326353).

\bibliographystyle{abbrv}
\bibliography{mybib}

\end{document}